\documentclass[11pt]{amsart}

\usepackage{amsmath}
\usepackage{amssymb}
\usepackage{mathrsfs}

\newtheorem{theorem}{Theorem}[section]
\newtheorem{claim}[theorem]{Claim}

\newtheorem{lemma}[theorem]{Lemma}

\theoremstyle{definition}
\newtheorem{definition}[theorem]{Definition}

\newtheorem{question}[theorem]{Question}

\theoremstyle{remark}

\newcount\skewfactor
\def\mathunderaccent#1#2 {\let\theaccent#1\skewfactor#2
\mathpalette\putaccentunder}
\def\putaccentunder#1#2{\oalign{$#1#2$\crcr\hidewidth
\vbox to.2ex{\hbox{$#1\skew\skewfactor\theaccent{}$}\vss}\hidewidth}}
\def\name{\mathunderaccent\tilde-3 }

\def\smallbox#1{\leavevmode\thinspace\hbox{\vrule\vtop{\vbox
   {\hrule\kern1pt\hbox{\vphantom{\tt/}\thinspace{\tt#1}\thinspace}}
   \kern1pt\hrule}\vrule}\thinspace}

\newcommand{\cf}{{\rm cf}}

\def\qedref#1{$\qed_{\reforiginal{#1}}$}

\title{Double weakness}
\author{Shimon Garti}
\address{Institute of Mathematics,
 The Hebrew University of Jerusalem,
 Jerusalem 91904, Israel}
\email{shimon.garty@mail.huji.ac.il}

\author{Saharon Shelah}
\address{Institute of Mathematics
 The Hebrew University of Jerusalem
 Jerusalem 91904, Israel
 and  Department of Mathematics
 Rutgers University
 New Brunswick, NJ 08854, USA}
\email{shelah@math.huji.ac.il}
\urladdr{http://www.math.rutgers.edu/\char`\~shelah}

\subjclass[2010]{03E35}
\keywords{Weak diamond, very weak diamond, weakly inaccessible, Cohen forcing, Radin forcing}
\thanks{This is publication 1111 of the second author. Both authors are grateful to the generous support of the European Research Council, grant no. 338821}

\begin{document}
\let\labeloriginal\label
\let\reforiginal\ref

\begin{abstract}
We prove that, consistently, there exists a weakly but not strongly inaccessible cardinal $\lambda$ for which the sequence $\langle 2^\theta:\theta<\lambda\rangle$ is not eventually constant and the weak diamond fails at $\lambda$.
We also prove that consistently diamond fails but a parametrized version of weak diamond holds at some strongly inaccessible $\lambda$.
\end{abstract}

\maketitle

\newpage

\section{Introduction}

The diamond on $\kappa$ is a prediction principle, discovered by Jensen, \cite{MR0309729}. Assuming that $\kappa=\cf(\kappa)>\omega$, the principle $\Diamond_\kappa$ says that there exists a sequence of sets $\langle A_\alpha:\alpha<\kappa\rangle$ so that $A_\alpha\subseteq\alpha$ for every $\alpha<\kappa$ and each subset $A$ of $\kappa$ is predicted by the diamond sequence at stationarily many points. That is, for every $A\subseteq\kappa$, the set $S_A=\{\alpha<\kappa:A_\alpha=A\cap\alpha\}$ is a stationary subset of $\kappa$.

There is an immediate implication of $\Diamond_\kappa$ on cardinal arithmetic below $\kappa$. If $\kappa=\theta^+$ then $\Diamond_\kappa$ implies $2^\theta=\theta^+$. 
If $\kappa=\theta^+>\aleph_1$ then the opposite implication holds as well, namely $2^\theta=\theta^+$ implies $\Diamond_\kappa$.
In fact, a more comprehensive statement holds and it applies to many stationary subsets of $\kappa$, see \cite{MR2596054}.
However, $2^{\aleph_0}=\aleph_1$ is consistent with $\neg\Diamond_{\aleph_1}$. This peculiarity motivated Devlin and Shelah, \cite{MR0469756}, to define a prediction principle sufficiently weak to follow from the continuum hypothesis:

\begin{definition}
\label{defweakdiamond} Weak diamond. \newline 
Let $\lambda$ be a regular and uncountable cardinal. \newline 
The weak diamond on $\lambda$, denoted $\Phi_\lambda$, is the assertion that for every coloring $c:{}^{<\lambda}2\rightarrow 2$ there exists a weak diamond function $g\in{}^\lambda 2$ such that for each element $f\in{}^\lambda 2$ the set $\{\alpha<\lambda:c(f\upharpoonright\alpha)=g(\alpha)\}$ is a stationary subset of $\lambda$.
\end{definition}

Concerning successor cardinals of the form $\lambda=\kappa^+$, the weak diamond can be characterized by a simple cardinal arithmetic statement. It is shown in several places that $\Phi_{\kappa^+}\Leftrightarrow 2^\kappa<2^{\kappa^+}$. However, there is no characterization for either diamond or weak diamond by pure cardinal arithmetic statements when $\lambda$ is regular and limit.

This fact motivated the definition of a yet weaker prediction principle called the very weak diamond and denoted $\Psi_\lambda$. The definition comes from \cite{MR3914938} and reads as follows:

\begin{definition}
\label{defvwd} Very weak diamond. \newline 
Let $\lambda$ be a regular and uncountable cardinal. \newline 
The very weak diamond on $\lambda$, denoted $\Psi_\lambda$, is the assertion that for every coloring $c:{}^{<\lambda}2\rightarrow 2$ there exists a \emph{very weak diamond} function $g\in{}^\lambda 2$ such that for each element $f\in{}^\lambda 2$ the set $\{\alpha<\lambda:c(f\upharpoonright\alpha)=g(\alpha)\}$ is an unbounded subset of $\lambda$.
\end{definition}

It has been proved in \cite{MR3914938} that $\Psi_\lambda$ is equivalent to $2^{<\lambda}<2^\lambda$ whenever $\lambda$ is regular and uncountable. In particular, if $\lambda$ is strongly inaccessible then $\Psi_\lambda$ follows. This is not the case for the elder brothers $\Diamond_\lambda$ and $\Phi_\lambda$. An unpublished result of Woodin says that $\neg\Diamond_\lambda$ is consistent for a strongly inaccessible cardinal $\lambda$, and his method has been applied in \cite{MR3914938} in order to prove a similar statement with respect to the weak diamond.

The purpose of this paper is two-fold. We shall begin with a survey of open problems, concerning the various principles mentioned above and their relationship with large cardinals. Then we shall address two of these problems.
Specifically, we separate diamond from a parametrized version of weak diamond at strongly inaccessibles and we settle the question of weakly inaccessibles by covering the last open case concerning these cardinals.

Our notation is mostly standard.
We shall use $\theta,\kappa,\lambda,\mu$ to denote cardinals and $\alpha,\beta,\gamma,\delta,\varepsilon,\zeta$ as well as $i,j$ to denote ordinals.
We denote ordinal multplication by $\alpha\cdot\beta$.
We adopt the Jerusalem forcing notation, that is $p\leq_{\mathbb{P}}q$ means that $p$ is weaker than $q$.
We shall use basic facts from pcf theory, and we suggest \cite{MR2768693}, \cite{MR1318912} and \cite{kojman} as a reference.
In the last section we force with Radin forcing, and we suggest \cite{MR2768695} as a reference.

\newpage 

\section{Open problems}

In what follows, a large cardinal is a regular and limit cardinal. This includes $\aleph_0$, but since all kinds of diamond principles fail at $\aleph_0$ we shall assume henceforth that $\lambda$ is uncountable. The behavior of $\Diamond_\lambda$ and $\Phi_\lambda$ is very interesting with respect to large cardinals. If $\lambda$ is large enough then $\Diamond_\lambda$ holds (a fortiori, $\Phi_\lambda$ and $\Psi_\lambda$).

For example, if $\lambda$ is measurable then $\Diamond_\lambda$ holds.
It has been proved by Kunen and Jensen \cite{jk} that if $\lambda$ is subtle then $\Diamond_\lambda$ holds, and one should bear in mind that subtlety is compatible with the constructible universe. However, if $\lambda$ is just strongly inaccessible or even strongly Mahlo, then $\neg\Diamond_\lambda$ is consistent (Woodin, \cite{woo}).

\begin{question}
\label{qdiamondwc} Is it consistent that $\lambda$ is weakly compact and $\neg\Diamond_\lambda$?
\end{question}

The method of Woodin is to use Radin forcing, starting from a measurable cardinal $\kappa$ with $2^\kappa>\kappa^+$. In the generic extension, $2^\alpha>\alpha^+$ at (almost) every $\alpha$ in the Radin club. In particular, \textsf{GCH} fails at unboundedly many places.

\begin{question}
\label{qdiamondgch} Is it consistent that $\lambda$ is strongly inaccessible, $\theta<\lambda\Rightarrow 2^\theta=\theta^+$ and $\neg\Diamond_\lambda$?
\end{question}

The negation of $\Phi_\lambda$ is a bit harder. Nevertheless, one can prove the consistency of $\neg\Phi_\lambda$ when $\lambda$ is strongly inaccessible or strongly Mahlo, as shown in \cite{MR3914938}. 

\begin{question}
\label{qwdwc} Is it consistent that $\lambda$ is weakly compact and $\neg\Phi_\lambda$?
\end{question}

By the same token:

\begin{question}
\label{qwdgch} Is it consistent that $\lambda$ is strongly inaccessible, $\theta<\lambda\Rightarrow 2^\theta=\theta^+$ and $\neg\Phi_\lambda$?
\end{question}

Another way to examine large cardinals around weak compactness is related to J\'onsson cardinals.
It is an open problem whether the first strongly inaccessible J\'onsson cardinal must be weakly compact.
On this ground, Dan Saattrup Nielsen, \cite{mof}, raised the following problem:

\begin{question}
\label{qjonsson} Is it consistent that $\lambda$ is strongly inaccessible and J\'onsson, yet $\neg\Diamond_\lambda$ or just $\neg\Phi_\lambda$?
\end{question}

Another question is connected with the fact that the failure of $\Diamond_\kappa$ in Woodin's model is witnessed by a set in the ground model.

\begin{question}
\label{qground} Is it consistent that $\lambda$ is strongly inaccessible, diamond fails at $\lambda$ but this cannot be witnessed by an element of $V$? That is, for some sequence of sets the only evidence for the failure of it to be a diamond sequence is an element of the generic extension.
\end{question}

An interesting issue arises with respect to weakly inaccessible cardinals. Suppose that $\lambda$ is such a cardinal. We distinguish three cases. In the first one $2^{<\lambda}=2^\lambda$, in which case $\Phi_\lambda$ fails (a fortiori, $\Diamond_\lambda$ fails). If $2^{<\lambda}=2^\kappa<2^\lambda$ for some $\kappa<\lambda$ then $\Phi_\lambda$ holds, this is the second case. The third case is when $2^{<\lambda}<2^\lambda$ but the sequence $\langle 2^\theta:\theta<\lambda\rangle$ is not eventually constant. Consistently, $\Phi_\lambda$ holds in this case (see \cite{MR3914938}).

\begin{question}
\label{qwdwinac} Is it consistent that $\lambda$ is weakly but not strongly inaccessible, the sequence $\langle 2^\theta:\theta<\lambda\rangle$ is not eventually constant, and $\neg\Phi_\lambda$?
\end{question}

Notice that such weakly inaccessible cardinals behave in some respect like strongly inaccessible cardinals (for which the sequence $\langle 2^\theta:\theta<\lambda\rangle$ is not eventually constant as well). However, it seems problematic to merge Radin forcing with blowing up the power set of small cardinals. 
Nevertheless, we shall give a positive answer to Question \ref{qwdwinac} in the next section, by proving that $\neg\Phi_\lambda$ is consistent in these cases.

We turn now to the interplay between $\Diamond_\lambda$ and $\Phi_\lambda$. For every successor cardinal $\lambda=\kappa^+$ it is easy to force $\neg\Diamond_\lambda\wedge\Phi_\lambda$ by forcing $\kappa^+<2^\kappa<2^{\kappa^+}$. We may ask:

\begin{question}
\label{qdwd} Is it consistent that $\lambda$ is strongly inaccessible, $\Phi_\lambda$ holds but $\Diamond_\lambda$ fails?
\end{question}

Golshani, \cite{go}, proved that the unpublished $\neg\Diamond_\lambda$ of Woodin as well as $\neg\Phi_\lambda$ of \cite{MR3914938} at strongly inaccessible cardinals can be pushed down to the first strongly inaccessible cardinal. The following is natural:

\begin{question}
\label{qfirst} Let $\lambda$ be the first strongly inaccessible cardinal. \newline 
Is it consistent that $\Phi_\lambda$ holds while $\Diamond_\lambda$ fails?
\end{question}

The statement $\neg\Diamond_\lambda$ where $\lambda$ is a smallish large cardinal has some consistency strength.
Jensen proved in \cite{je91} that if $\lambda$ is Mahlo and $S=S^\lambda_{\omega_1}$ then $\neg\Diamond_S$ implies that $0^\sharp$ exists.
Zeman, \cite{MR1812181}, employed these ideas to show that if $\lambda$ is Mahlo, $\theta=\cf(\theta)<\lambda$ and $S=S^\lambda_\theta$ then $\neg\Diamond_S$ implies the existence of many measurable cardinals $\kappa<\lambda$ with $o(\kappa)\geq\theta$ in the core model $K$.
It seems, however, that the methods of Jensen and Zeman apply to Mahlo cardinals.

\begin{question}
\label{qinacstrength} What is the consistency strength of $\neg\Diamond_\lambda$ where $\lambda$ is strongly inaccessible and not Mahlo?
\end{question}

A comparison between $\Phi_\lambda$ and $\Diamond_\lambda$ from this point of view is expressed by the following:

\begin{question}
\label{qconstr} What is the consistency strength of $\neg\Phi_\lambda$ where $\lambda$ is strongly inaccessible? 
In particular, is it strictly stronger than $\neg\Diamond_\lambda$ where $\lambda$ is strongly inaccessible?
\end{question}

If one wishes to focus on specific stationary sets, then interesting results are known above weak compactness. For a strongly Mahlo cardinal $\kappa$ let $S_M(\kappa)$ be the set $\{\theta<\kappa:\theta=\cf(\theta)\}$. Woodin proved the consistency of $\neg\Diamond_{S_M(\kappa)}$ when $\kappa$ is weakly compact. Hauser, \cite{MR1164732}, stressed the method of Woodin toward indescribable cardinals, and D\v zamonja-Hamkins \cite{MR2279655} contains the same consistency negative result at unfoldable cardinals.

\begin{question}
\label{qwdlarge} Let $\kappa$ be weakly compact. \newline 
Is it consistent that $\neg\Phi_{S_M(\kappa)}$?
\end{question}

\newpage

\section{Weak diamond at weakly inaccessibles}

In this section we address Question \ref{qwdwinac} by proving the consistency of the negation of the weak diamond at a weakly inaccessible cardinal $\lambda$ which behaves like a strongly inaccessible cardinal in the sense that $\langle 2^\theta:\theta<\lambda\rangle$ is not eventually constant. Notice that in such a case we have $2^{<\lambda}<2^\lambda$, hence the very weak diamond $\Psi_\lambda$ holds. 
The fact that $\Psi_\lambda$ is strictly weaker than $\Phi_\lambda$ has been already proved in \cite{MR3914938} in a different way, thus we have another example to the fact that $\Psi_\lambda$ is strictly weaker than $\Phi_\lambda$.

Ahead of proving the main result we phrase the following.
The proof is straightforward, using basic facts from pcf theory.
Probably the statement holds in \textsf{ZFC}, but since we use the claim in order \emph{to force} $\neg\Phi_\lambda$, this is unimportant.

\begin{claim}
\label{clmpcf} Assume that $\langle\mu_i:i\leq\lambda\rangle$ is an increasing continuous sequence of cardinals, $\lambda<\mu_0$ and $2^{\mu_i}=\mu_i^+$ for every limit ordinal $i\leq\lambda$. \newline 
Then there is a $\mu_0^+$-complete forcing notion $\mathbb{Q}$ which neither collapses cardinals nor changes cofinalities and if $G\subseteq\mathbb{Q}$ is $V$-generic then there are $\bar{f},\bar{d}\in V[G]$ such that:
\begin{enumerate}
\item [$(a)$] $\bar{f}=\langle\bar{f}_i:i\leq\lambda, i\ \text{is a limit ordinal}\rangle$.
\item [$(b)$] $\bar{f}_i=\langle f_{i,\alpha}:\alpha\in\mu_i^+\rangle$.
\item [$(c)$] $f_{i,\alpha}\in\prod_{j<i}\mu_j^+$ for every $\alpha\in\mu_i^+$.
\item [$(d)$] $\bar{f}_i$ is $J^{\rm bd}_i$-increasing with $\alpha$, and cofinal in the product $\prod_{j<i}\mu_j^+$. 
\item [$(e)$] $\bar{d}=\langle d_i:i<\lambda\rangle$.
\item [$(f)$] $d_i$ is a function from $\mu_i^+$ into $\mu_i^+$, for every $i\in\lambda$.
\item [$(g)$] For every limit ordinal $i\leq\lambda$, every $\eta\in\prod_{j<i}\mu_j^+$ and arbitrarily large $\alpha\in\mu_i^+$, the following statement holds: for every limit ordinal $j<i$ we have $d_j(\beta)=\eta(i)$, where $\beta$ is the unique ordinal in $\mu_j^+$ such that $f_{\lambda,\alpha}\upharpoonright i=f_{j,\beta}$.
\end{enumerate}
\end{claim}

\hfill \qedref{clmpcf}

Now we can state the main theorem:

\begin{theorem}
\label{thmwdwc} Let $\lambda$ be strongly inaccessible. \newline 
There is a forcing notion $\mathbb{P}$ which forces $\lambda$ to be weakly but not strongly inaccessible, the sequence $\langle 2^\theta:\theta<\lambda\rangle$ is not eventually constant and $\neg\Phi_\lambda$.
\end{theorem}

\par\noindent\emph{Proof}. \newline 
We begin with the following assumptions (in the ground model) on the strongly inaccessible cardinal $\lambda$:
\begin{itemize}
\item $\langle\lambda_i:i<\lambda\rangle$ is an increasing and continuous sequence of cardinals, $\lambda=\bigcup\{\lambda_i:i<\lambda\}$ and $\lambda_0=0$.
\item If $i<\lambda$ is a limit ordinal then $2^{\lambda_i}=\lambda_i^+=\lambda_{i+1}$, and $2^\lambda=\lambda^+$.
\item If $i<\lambda$ is a successor ordinal then $\lambda_i=\lambda_i^{<\lambda_i}$ and $\alpha<\lambda_{i+1} \Rightarrow |\alpha|^{\lambda_i}<\lambda_{i+1}$.
\item $\langle\mu_i:i\leq\lambda\rangle$ is an increasing and continuous sequence of cardinals, $\mu=\mu_\lambda$ and $\mu_0=0<\lambda<\mu_1$.
\item If $i<\lambda$ is a limit ordinal then $\mu_i^\lambda=\mu_i^+=\mu_{i+1}$, and $2^{\mu_\lambda}=\mu^+_\lambda$.
\item If $i<\lambda$ is a successor ordinal then $\mu_i$ is regular and $\alpha<\mu_{i+1} \Rightarrow |\alpha|^{\lambda_i}<\mu_{i+1}$.
\end{itemize}

The idea is to add for an unbounded set of $\lambda_i$-s many Cohen subsets (\emph{many} is determined here by the $\mu_i$-s), thus making $\lambda$ weakly but not strongly inaccessible cardinal. Likewise, the sequence of $2^{\lambda_i}$-s will not be eventually constant. Of course, the main task is to make sure that $\neg\Phi_\lambda$ holds in the generic extension.

Firstly, we force with $\mathbb{Q}$ from Claim \ref{clmpcf}, thus $\bar{f}$ and $\bar{d}$ are as described there.
Secondly, we define a forcing notion $\mathbb{P}$ which will force the desired properties of $\lambda$ and concomitantly $\neg\Phi_\lambda$.
We define the basic component, the forcing notion $\mathbb{P}_i$ for every $i<\lambda$, as ${\rm Add}(\lambda_{i+2},\mu_{i+2})$. Explicitly, a condition $p\in\mathbb{P}_i$ is a partial function from $\mu_{i+2}$ into $\lambda_{i+2}$ such that $|{\rm dom}(p)|<\lambda_{i+2}$. If $p,q\in\mathbb{P}_i$ then $p\leq_{\mathbb{P}_i} q$ iff $p\subseteq q$. For every $j<\lambda$ let $\mathbb{P}_{<j}$ be the product (not an iteration) $\prod\limits_{i<j} \mathbb{P}_i$ with full support. Likewise, $\mathbb{P}=\mathbb{P}_{<\lambda}$. The forcing notion $\mathbb{P}$ will satisfy the statements of the theorem.

Suppose that $i<\lambda$ and let $G(i)$ be a generic subset of $\mathbb{P}_i$. By the definition of $\mathbb{P}_i$ we have a name $\name{\eta}_i^{G(i)} = \name{\eta}_i\in {}^{\mu_{i+2}} \lambda_{i+2}$. So $\name{\eta}_i$ is a name for a function into $\lambda_{i+2}$, but occasionally we would like to restrict the range of the function.
Hence we define the following version. For every $i<\lambda$ let $\name{\eta}^1_i\in{}^{\mu_{i+2}} 2$ be defined by $\name{\eta}^1_i(\beta) = \name{\eta}_i(\beta)$ if $\name{\eta}_i(\beta)<2$ and zero otherwise.

We point out some basic properties of the above forcing notions. Firstly we observe that $\mathbb{P}_{<i}$ is $\lambda_{i+1}^+$-cc and $\mathbb{P}_i$ is $\lambda_{i+2}^+$-cc for every $i<\lambda$. This follows from a usual Delta-system argument, by recalling that $\lambda_{i+1}=\lambda_{i+1}^{<\lambda_{i+1}}$. 
Secondly, for every $i<\lambda$ we can decompose the product $\mathbb{P}$ into two parts, namely $\mathbb{P}_{<i}$ and $\mathbb{P}_{[i,\lambda)} = \prod\limits_{j\in[i,\lambda)} \mathbb{P}_j$. 
It is easily verified that $\mathbb{P} \cong \mathbb{P}_{<i}\times \prod\limits_{j\in[i,\lambda)} \mathbb{P}_j$ for every $i<\lambda$. 
Indeed, for any condition $p\in\mathbb{P}$ let $q=p\upharpoonright i$ and $r=p-p\upharpoonright i$, so $q\in\mathbb{P}_{<i}$ and $r\in\mathbb{P}_{[i,\lambda)}$.
Now the mapping $h:\mathbb{P} \rightarrow \mathbb{P}_{<i}\times \prod\limits_{j\in[i,\lambda)} \mathbb{P}_j$ given by $h(p)=(q,r)$ is an isomorphism of forcing notions.

Notice, further, that $\mathbb{P}_{[i,\lambda)}$ is $\lambda_{i+2}$-complete for every $i<\lambda$. This follows from the regularity of $\lambda_{i+2}$ and the fact that $\mathbb{P}_{[i,\lambda)}$ is a product, hence inherits the completeness of its components. 
Finally, observe that $\mathbb{P}$ is $\lambda^{++}$-cc, so every cardinal above $\lambda^{++}$ is preserved by $\mathbb{P}$. 
It is also evident that $\lambda^+$ is preserved upon forcing with $\mathbb{P}$.
We conclude that forcing with $\mathbb{P}$ preserves cardinals, no cofinality is changed and if $G\subseteq\mathbb{P}$ is generic then $V[G]\models \forall i<\lambda, 2^{\lambda_{i+2}}=\mu_{i+2} \wedge 2^\lambda=\mu_\lambda^+=\mu^+$.
Therefore, $\lambda$ is weakly but not strongly inaccessible, and $\langle 2^\theta:\theta\in\lambda\rangle$ is not eventually constant in the generic extension $V[G]$.
It remains to show that $V[G]\models\neg\Phi_\lambda$.

Assume that $p\in\mathbb{P}$ satisfies $p\Vdash_{\mathbb{P}} \name{\tau}: \check{\lambda}\rightarrow Ord$. Observe that one can find an extension $p\leq q\in\mathbb{P}$ so that:
\begin{enumerate}
\item [$(a)$] For every $i<\lambda$ there is a name $\name{\tau}_i\in \mathbb{P}_{<i}$ such that $q\Vdash\name{\tau}\upharpoonright\lambda_{i+1} = \name{\tau}_i \wedge q\upharpoonright i\Vdash_{\mathbb{P}_{<i}} \name{\tau}_i: \check{\lambda}_i\rightarrow Ord$.
\item [$(b)$] If $i<\lambda$ is a limit ordinal then $\name{\tau}_i$ is a $\mathbb{P}_{<i}$-name which is determined by $\lambda_{i+1}$ maximal antichains in $\mathbb{P}_{<i}$ of the form $\bar{p}_{i,\varepsilon} = \langle p_{i,\varepsilon,\alpha}: \alpha<\alpha_{i,\varepsilon}\leq\lambda_i^+ \rangle$ for each $\varepsilon<\lambda_i$.
\end{enumerate}

For every limit ordinal $i<\lambda$, if $\mathcal{A}\subseteq\mathbb{P}_{<i}$ is a maximal antichain then $|\mathcal{A}|\leq\lambda_{i+1}$.
Hence the number of nice $\mathbb{P}_{<i}$-names of an element of $\{0,1\}$ is bounded by the number of antichains of $\mathbb{P}_{<i}$ which is at most $\mu_i^+\cdot 2^{\lambda_{i+1}} = \mu_i^+=\mu_{i+1}$.
So for every limit ordinal $i<\lambda$, let $\bar{\sigma}_i=\langle\name{\sigma}_{i\alpha}: \alpha\in\mu_{i+1}\rangle$ list all the nice $\mathbb{P}_{<i}$-names of an element of $\{0,1\}$.

The main idea of the proof is to enumerate all the names for possible weak diamond sequences, and then to eliminate each candidate by forcing the opposite value at an end segment, using the conditions of $\mathbb{P}$. For this end, we shall define the concept of comprehensive sequences as follows.
For every limit ordinal $i<\lambda$ we choose a sequence of names $\name{\bar{\tau}}_i=\langle \name{\tau}_{i,\beta}:\beta<\mu_i^+\rangle$, with the following properties:
\begin{enumerate}
\item [$(\aleph)$] Each $\name{\tau}_{i,\beta}$ is a nice name in $\mathbb{P}_{<i}$ of an element from ${}^{\lambda_i}2$.
\item [$(\beth)$] For every $\varepsilon<\lambda_i$ there is a maximal antichain $\bar{q}_{i,\varepsilon} = \langle q_{i,\varepsilon,\beta}: \beta<\mu_i^+ \rangle$ in the forcing $\mathbb{P}_{<i}$ (the length, without loss of generality, is exactly $\mu_i^+$), and a sequence of values $\bar{t}_{i,\varepsilon} = \langle t_{i,\varepsilon,\beta}:\beta<\mu_i^+ \rangle$ so that $q_{i,\varepsilon,\beta}\Vdash \name{\tau}_{i,\varepsilon}(\beta) = t_{i,\varepsilon,\beta}$.
\item [$(\gimel)$] If $\beta_0<\beta_1<\mu_i^+$ then $\Vdash_{\mathbb{P}_{<i}} \name{\tau}_{i,\beta_0}\neq \name{\tau}_{i,\beta_1}$.
\item [$(\daleth)$] It is forced by $\mathbb{P}_{<i}$ that every element of ${}^{\lambda_i}2$ appears in the list of sequences, i.e. if $p_0\Vdash_{\mathbb{P}_{<i}} \name{\sigma}\in {}^{\lambda_i}2$ then there exist an ordinal $\beta<\mu_i^+$ and a condition $p_1\in \mathbb{P}_{<i}$ so that $p_0\leq p_1$ and $p_1\Vdash \name{\sigma}=\name{\tau}_{i,\beta}$.
\end{enumerate}

Given a name of such a sequence and a condition which forces this fact, we wish to make sure that there is a large set of free points in which we can extend our condition and force any desired value.

For every limit ordinal $i\leq\lambda$ and every function $f\in\bar{f}_i$ we define a $\mathbb{P}_{<i}$-name $\name{\nu}_f$ of a mapping from $\mu_i$ into $2$ by the following procedure. If $j<i$ and $\alpha\in[\mu_j,\mu_{j+1})$ then $\name{\nu}_f(\alpha)=\name{\eta}^1_j(f(j+1)+\alpha)$. 

For proving $\neg\Phi_\lambda$ we have to define a coloring $\name{c}:{}^{<\lambda}2\rightarrow 2$ which is not predicted by any candidate of a weak diamond sequence.
We define $\name{c}$ in a piecewise manner, so for every limit ordinal $i<\lambda$ we define a $\mathbb{P}_{<i}$-name $\name{c}_i$ of a mapping from $({}^{\lambda_i}2)^{V[\mathbb{P}_{<i}]}$ into $2$, as follows.

We shall define $\name{c}_i(\name{\tau})$ as a $\mathbb{P}$-name of an element in $\{0,1\}$ for every limit ordinal $i<\lambda$ and every $\mathbb{P}_{<i}$-name $\name{\tau}$ of an element in ${}^{\lambda_i}2$.
Define $\name{c}_i(\name{\tau})=1-\name{\sigma}_{i\varepsilon}$ if $\name{\tau}=\name{\nu}_{f_{i,\varepsilon}}$ for some $\varepsilon<\mu_i^+$, and zero otherwise.

Let $\name{c}=\bigcup\{\name{c}_i:i<\lambda, i\ \text{is a limit ordinal}\}$.
Let $E=\{\lambda_i:i<\lambda, i$ is a limit ordinal$\}$.
It is enough to prove that there is no guessing function $h\in{}^\lambda 2$ in $V^{\mathbb{P}}$ witnessing $\Phi_E$ with respect to $\name{c}$.
Assume towards contradiction that there are $p_\ast\in\mathbb{P}$ and $\name{h}\in{}^E 2$ so that $p_\ast$ forces $\name{h}$ to be a weak diamond guessing function for $\name{c}$.
Let $\name{\sigma}_i^*$ be a $\mathbb{P}_{<i}$-name of an element of $\{0,1\}$ (for every limit ordinal $i<\lambda$) such that $p_\ast(i)\Vdash\name{h}(\lambda_i)=\name{\sigma}_i^*$.
We may assume, without loss of generality, that this happens for every $i$.

By the choice of the sequence $\bar{\sigma}_i$ we can find, for every limit $i<\lambda$, an ordinal $\beta_i\in\mu_i^+$ for which $\name{\sigma}_i^*=\name{\sigma}_{i\beta_i}$.
For every limit ordinal $i<\lambda$ define $A_i=\{\beta\in\mu_{i+1}:\name{\sigma}_{i\beta}=\name{\sigma}_i^*\}$.
Notice that $A_i$ is an unbounded subset of $\mu_{i+1}$.
By the choice of the scale $\langle f_\zeta:\zeta\in\mu^+\rangle$ and the objects given in Claim \ref{clmpcf} we can see that for every sufficiently large limit ordinal $i<\lambda$ there are $\zeta_i,\varepsilon_i$ such that $f_{\lambda,\zeta_i}\upharpoonright i=f_{i,\varepsilon_i}$.
By the scale properties there are $\zeta\in\mu^+$ and $i(\ast)\in\lambda$ such that if $i\in[i(\ast),\lambda)$ then $\sup({\rm dom}p_\ast(i))<f_\zeta(i)$.

Consider the $\mathbb{P}$-name $\name{\nu}_{f_\zeta}$ defined above.
It is a name of an element of ${}^\lambda 2$.
We claim that there is a condition $q\in\mathbb{P}$ so that $p_\ast\leq q$ and if $i\in[i(\ast),\lambda)$ is a limit ordinal then $q\Vdash\name{c}_i(\name{\nu}_{f_\zeta}\upharpoonright\lambda_i)\neq \name{\sigma}_i^*$.
To show that there is such a condition $q$ recall the definition of $\name{c}_i$ and find for each limit ordinal $i\in\lambda$ an ordinal $\varepsilon_i\in\mu_i^+$ such that $f_{\lambda,\zeta}\upharpoonright\lambda_i = f_{i,\varepsilon_i}$.
Now extend $p_\ast$ to a condition $q$ which satisfies $q\Vdash_{\mathbb{P}} \name{\sigma}_{i,d_i(\varepsilon_i)}=1-\name{\sigma}_{i,\beta_i}$ for every sufficiently large $i\in\lambda$.

It follows that $q\Vdash\name{h}(\lambda_i)=\name{\sigma}_i^*=\name{\sigma}_{i,\beta_i} \neq 1-\name{\sigma}_{i,\beta_i}=\name{\sigma}_{i,d_i(\varepsilon_i)}= \name{c}_i(\name{\nu}_{f_\zeta}\upharpoonright\lambda_i)$ for an end-segment of elements of $E$, so we are done.

\hfill \qedref{thmwdwc}

An iteration of Cohen forcing has been used in \cite{MR3914938} in order to prove $\Phi_\lambda$ where $\lambda$ is weakly but not strongly inaccessible and $\langle 2^\theta:\theta\in\lambda\rangle$ is not eventually constant.
Hence Cohen forcing may produce the positive statement $\Phi_\lambda$ as well as the negative statement $\neg\Phi_\lambda$.
The main point is the support of the iteration (or the product).
For proving the positive direction one uses Easton support, while here we have used full support.

We conclude with a short remark about consistency strength.
Woodin's method to force $\neg\Diamond_\lambda$ at small large cardinals employs Radin forcing and hence requires some hyper-measurability assumptions.
Zeman proved in \cite{MR1812181} that such assumptions are necessary.
However, for the main result of this section one can begin with the constructible universe and with one strongly inaccessible cardinal.

\newpage 

\section{At strongly inaccessibles}

In this section we prove a partial result concerning Question \ref{qdwd}.
We shall show that a parametrized version of weak diamond, which is weaker than the traditional weak diamond, holds in Woodin's model of $\neg\Diamond_\lambda$ where $\lambda$ is strongly inaccessible.

\begin{definition}
\label{defparametrized} Parametrized weak diamond. \newline 
Let $\theta,\kappa$ be cardinals, where $\kappa=\cf(\kappa)>\aleph_0$. \newline 
We shall say that $\theta-\Phi_\kappa$ holds iff for every coloring $c:{}^{<\kappa}2\rightarrow 2$ there is a family $E=\{g_\tau:\tau\in\theta\}\subseteq{}^\kappa 2$ such that for every $f\in{}^\kappa 2$ one can find $\tau\in\theta$ for which $\{\alpha\in\kappa:c(f\upharpoonright\alpha)=g_\tau(\alpha)\}$ is a stationary subset of $\kappa$.
\end{definition}

The usual weak diamond $\Phi_\kappa$ is $1-\Phi_\kappa$ in the above terminology.
The interesting instances of $\theta-\Phi_\kappa$ are only when $\theta<2^\kappa$.
We shall prove that in Woodin's model for $\neg\Diamond_\lambda$ at a strongly inaccessible $\lambda$ one has $\lambda^+-\Phi_\lambda$ where $\lambda^+<2^\lambda$.

Woodin's method is based on Radin forcing, which appeared in \cite{MR670992}.
There are several ways to describe this forcing notion, and we take Mitchell's approach from \cite{MR673794}.
We shall use the notational conventions of \cite{MR2768695}.
However, we sketch, briefly, some basic facts and fix our notation.

Let $\kappa$ be a measurable cardinal with $o(\kappa)\geq\kappa^+$.
We shall force with a measure sequence $\bar{U}$ which comes from an elementary embedding, so we begin with $\jmath:V\rightarrow M$ so that $\kappa={\rm crit}(\jmath)$ and the measures are derived from $\jmath$ by the following procedure.
If $\tau\in\ell(\bar{U})$ then $U_\tau=\{A\subseteq V_\kappa:\bar{U}\upharpoonright\tau\in\jmath(A)\}$.
We denote the sequence by $\langle\kappa\rangle^\frown\bar{U}$, so actually the length of the sequence is $1+\ell(\bar{U})$.

Every sequence of this form will be called \emph{a measure sequence}, and if $\kappa$ is measurable then $MS(\kappa)=\{\bar{\nu}:\bar{\nu}$ is a measure sequence on $\nu, \nu<\kappa\}$.
Notice that if $\bar{\nu}\in MS(\kappa)$ and $\bar{\nu}=
\langle\nu\rangle^\frown\langle\nu_\beta:\beta\in\ell(\bar{\nu})\rangle$ then each $\nu_\beta$ is a measure over $V_\nu$.
We let $\nu=\kappa(\bar{\nu})$, the measurable cardinal associated with $\bar{\nu}$.
We also denote $\bigcap\{\nu_\beta:\beta\in\ell(\bar{\nu})\}$ by $\bigcap\bar{\nu}$.
Following Gitik's conventions, $MS(\kappa)=A^0$ and for each $n\in\omega$ we let $A^{n+1}=\{\bar{\nu}\in A^n: A^n\cap V_{\kappa(\bar{\nu})}\in\bigcap\bar{\nu}\}$.
The set $\bar{A}=\bigcap\{A^n:n\in\omega\}$ is therefore an element of $\bigcap\bar{U}$ and we always assume that our large sets are subsets of $\bar{A}$.

\begin{definition}
\label{defradin} Radin forcing. \newline 
Let $\bar{U}$ be a measure sequence over $\kappa$.
A condition $p\in\mathbb{R}(\bar{U})$ is a finite sequence of the form $p=\langle d_i:i\leq\ell\rangle$, such that:
\begin{enumerate}
\item [$(\aleph)$] If $i\in\ell$ then either $d_i=\langle\kappa_i\rangle$ where $\kappa_i\in\kappa$, or $d_i=\langle\bar{\nu}_i,a_i\rangle$ where $\bar{\nu}_i$ is a measure sequence over $k_i=\kappa(\bar{\nu}_i)\in\kappa$ and $a_i\in\bigcap\bar{\nu}_i$.
\item [$(\beth)$] The sequence $\langle\kappa_i:i\in\ell\rangle$ is strictly increasing.
\item [$(\gimel)$] $d_\ell=\langle\bar{U},A\rangle$, where $A\in\bigcap\bar{U}$ and $A\subseteq\bar{A}$.
\item [$(\daleth)$] If $i<j\leq\ell$ and $d_j=\langle\bar{\nu}_j,a_j\rangle$ then $a_j\cap V_{\kappa_i+1}=\varnothing$.
\end{enumerate}
\end{definition}

If $p\in\mathbb{R}(\bar{U})$ then we usually introduce $p$ as $p_0^\frown\langle\bar{U},A\rangle$ where $p_0=\langle d_i:i\in\ell\rangle$ and $d_\ell=\langle\bar{U},A\rangle$.
We use the notation $\kappa(d_i),a(d_i)$ and $\bar{\nu}(d_i)$ to denote the pertinent objects mentioned in $d_i$.

\begin{definition}
\label{defradinorder} The pure order and the forcing order. \newline 
Assume that $p,q\in\mathbb{R}(\bar{U})$.
\begin{enumerate}
\item [$(\aleph)$] $p\leq^* q$ iff $p=\langle d_i:i\leq\ell\rangle, q=\langle e_i:i\leq\ell\rangle, \kappa(d_i)=\kappa(e_i)$ for every $i\leq\ell$ and $a(e_i)\subseteq a(d_i)$ whenever exists.
\item [$(\beth)$] $q$ is a one-point extension of $p$ iff for some $j\leq\ell$ and $\bar{\nu}\in a(d_j)$ we have $\kappa(\bar{\nu})>\kappa(d_{j-1})$ and either $\bar{\nu}=\langle\alpha\rangle$ and then $q=\langle d_i:i<j\rangle^\frown\langle\bar{\nu}\rangle^\frown\langle d_i:i\leq j\rangle$ or $\bar{\nu}$ contains measures on $V_{\kappa(\bar{\nu})}$ and then $q=\langle d_i:i<j\rangle^\frown\langle(\bar{\nu},a(d_j)\cap V_{\kappa(\bar{\nu})})\rangle^\frown\langle d_i:i\leq j\rangle$.
\item [$(\gimel)$] $p\leq q$ iff $q$ is obtained from $p$ by a finite process of pure extensions and one-point extensions.
\end{enumerate}
\end{definition}

We shall use the notation $q=p+\langle\bar{\nu}\rangle$ when $q$ is a one-point extension of $p$.
If $G\subseteq\mathbb{R}(\bar{U})$ is $V$-generic then $MS_G$ is the set of all measure sequences mentioned in $p$ for some $p\in G$, apart from the sequence $\bar{U}$ itself.
The set $C_G=\{\kappa(\bar{\nu}):\bar{\nu}\in MS_G\}$ is closed and unbounded in $\kappa$ in $V[G]$, and will be called \emph{the Radin club}.

The following lemma says that under the appropriate assumption one can take care of many colorings simultaneously.

\begin{lemma}
\label{lemsimul} Let $\kappa$ be a measurable cardinal. \newline 
Suppose that $c_\gamma:{}^{<\kappa}2\rightarrow 2$ is given for every $\gamma\in\kappa$.
Then there exists a weak diamond function $g\in{}^\kappa 2$ such that for every $f\in{}^\kappa 2$ and every $\gamma\in\kappa$ the set $\{\alpha\in\kappa:c_\gamma(f\upharpoonright\alpha)=g(\alpha)\}$ is a stationary subset of $\kappa$.
\end{lemma}

\par\noindent\emph{Proof}. \newline 
Let $\{S_\gamma:\gamma\in\kappa\}$ be a partition of $\kappa$ into $\kappa$-many disjoint stationary sets so that $\Phi_{S_\gamma}$ holds at every $S_\gamma$ (this is always possible if $\kappa$ is measurable).
For every $\gamma\in\kappa$ let $g_\gamma\in{}^{S_\gamma}2$ be a weak diamond function for the coloring $c_\gamma$.
Define $g=\bigcup_{\gamma\in\kappa}g_\gamma$, and notice that $g\in{}^\kappa 2$ is well defined since the domains of the $g_\gamma$s are mutually disjoint.
One can verify now that $g$ satisfies the statement of the lemma.

\hfill \qedref{lemsimul}

Now we can state the main result of this section:

\begin{theorem}
\label{thmq1.7} Separating diamond from weak diamond. \newline 
If there is a measurable cardinal $\kappa$ with $o(\kappa)\geq\kappa^+$ then one can force $\kappa$ to be strongly inaccessible, $\kappa^+<2^\kappa$ and $\kappa^+-\Phi_\kappa\wedge\neg\Diamond_\kappa$ hold simultaneously.
\end{theorem}

\par\noindent\emph{Proof}. \newline 
Let $\jmath:V\rightarrow M$ be an elementary embedding with $\kappa={\rm crit}(\jmath)$, and let $\bar{U}$ be the derived measure sequence.
We assume that $\ell(\bar{U})=\kappa^+$ and we force with $\mathbb{R}(\bar{U})$.
We assume, further, that $\kappa^+<2^\kappa$ in the ground model.
Let $G\subseteq\mathbb{R}(\bar{U})$ be $V$-generic.
We claim that $V[G]\models\Phi_\kappa\wedge\neg\Diamond_\kappa$.
Since $\kappa$ remains strongly inaccessible in $V[G]$, this will accomplish the proof.

The fact that $V[G]\models\neg\Diamond_\kappa$ follows from the assumption $\kappa^+<2^\kappa$ in $V$, and this is Woodin's argument.
We unfold the details since this will be the prototype for deriving $\kappa^+-\Phi_\kappa$.

Let $\name{s}=(\name{s}_\alpha:\alpha\in\kappa)$ be an $\mathbb{R}(\bar{U})$-name of a candidate of a diamond sequence.
Let $p\in\mathbb{R}(\bar{U})$ be a condition which forces $\name{s}_\alpha\subseteq\check{\alpha}$ for every $\alpha\in\kappa$.
We must find a condition $q\geq p$ which forces that $\name{s}$ fails to predict all the subsets of $\kappa$.
Specifically, we will choose a set $x\subseteq\kappa$ and a club $C\subseteq\kappa$ so that $q$ forces $\check{x}\cap\alpha\neq\name{s}_\alpha$ whenever $\alpha\in C$.
Remark that $x\in V$, so Woodin's argument gives more than the failure of $\Diamond_\kappa$ by showing that the failure can be witnessed by an element of the ground model.

Firstly, let $p=p_0^\frown(\bar{U},A^p)$, so $p\leq p+\langle\bar{\nu}\rangle$ for every $\bar{\nu}\in A^p$.
It will be useful to replace the name $\name{s}_\alpha$ (where $\alpha$ is $\kappa(\bar{\nu})$ for some $\bar{\nu}\in A^p$) by a name $\name{t}_\alpha$ from a smaller part of $\mathbb{R}(\bar{U})$.
For this end, we choose for each $\bar{\nu}\in A^p$ an $\mathbb{R}(\bar{\nu})$-name $\name{t}_{\kappa(\bar{\nu})}$ and a condition $(\bar{U},A_{\bar{\nu}})$ such that $(\bar{U},A_{\bar{\nu}})\Vdash \name{s}_{\kappa(\bar{\nu})}=\name{t}_{\kappa(\bar{\nu})}$.
We define $A^*=\Delta\{A_{\bar{\nu}}:\bar{\nu}\in A^p\}$, so $A^*\in\bigcap\bar{U}$.
Let $r=p_0^\frown(\bar{U},A^*)$, so $p\leq^* r$ and $r+\langle\bar{\nu}\rangle\Vdash \name{s}_{\kappa(\bar{\nu})}=\name{t}_{\kappa(\bar{\nu})}$ whenever $\bar{\nu}\in A^*$.
Our desired condition $q$ will be a pure extension of $r$.

For defining $q$ let $f:A^*\rightarrow V$ be the function $f(\bar{\nu})=\name{t}_{\kappa(\bar{\nu})}$.
Applying $\jmath$ to the formula $r+\langle\bar{\nu}\rangle\Vdash f(\bar{\nu})=\name{t}_{\kappa(\bar{\nu})}$ we see that:
$$
\forall\tau\in\kappa^+, \jmath(r)+\langle\bar{U}\upharpoonright\tau\rangle \Vdash\jmath(\name{s})_\kappa=\jmath f(\bar{U}\upharpoonright\tau).
$$
Now $\jmath f(\bar{U}\upharpoonright\tau)$ is a name of a subset of $\kappa$ in $M$.
Since we are assuming that $\kappa^+<2^\kappa$, there is some $x\subseteq\kappa$ so that $\jmath(r)+\langle\bar{U}\upharpoonright\tau\rangle \Vdash\forall\tau\in\kappa^+, \jmath f(\bar{U}\upharpoonright\tau)\neq\check{x}$.
Indeed, for every $\tau\in\kappa^+$ there is an antichain $\mathcal{A}_\tau$ of size at most $\kappa$ from which the generic set $G$ chooses the interpretation of $\jmath(\name{s})_\kappa$.
The size of the set of all the elements of $\mathcal{A}_\tau$ for every $\tau\in\kappa^+$ is $\kappa^+$ which is less than $2^\kappa$, so one can choose $x\subseteq\kappa$ as required.

By elementarity, there is a condition $q=p_0^\frown(\bar{U},A^q)$ so that $r\leq^* q$ and:
$$
\forall\bar{\nu}\in A^q, q+\langle\bar{\nu}\rangle\Vdash \name{s}_{\kappa(\bar{\nu})}\neq\check{x}\cap\kappa(\bar{\nu}).
$$
Since $\ell(\bar{U})=\kappa^+$ we see that $\kappa$ remains strongly inaccessible in $V[G]$, and $C=\{\kappa(\bar{\nu}):\bar{\nu}\in A^q\}\subseteq C_G$ is a club subset of $\kappa$.
The fact that $q\Vdash\forall\alpha\in C,\name{s}_\alpha\neq\check{x}\cap\alpha$ means that $\name{s}$ fails to name a diamond sequence, and since $p,\name{s}$ were arbitrary we conclude that $V[G]\models\neg\Diamond_\kappa$.

Moving to the weak diamond, assume that $\name{c}:{}^{<\kappa}2\rightarrow 2$ is a name of a coloring where $\name{c}=(\name{c}_\alpha:\alpha\in\kappa)$ and $\name{c}_\alpha:{}^{<\alpha}2\rightarrow 2$ for every $\alpha\in\kappa$.
Fix a condition $p\in\mathbb{R}(\bar{U})$ which forces this fact, and express $p$ as $p_0^\frown(\bar{U},A^p)$.
We are trying to find $q\geq p$ and a set of functions $E=\{g_\gamma:\gamma\in\kappa^+\}\subseteq{}^\kappa 2$ such that $q$ forces $E$ to be a witness to $\kappa^+-\Phi_\kappa$ for $\name{c}$.

As in the first part of the proof we shrink $A^p$ by the following procedure.
For each $\bar{\nu}\in A^p$ we choose an $\mathbb{R}(\bar{\nu})$-name $\name{d}_{\kappa(\bar{\nu})}$ and a set $A_{\bar{\nu}}\subseteq A^p$ such that $(\bar{U},A_{\bar{\nu}})\Vdash \name{c}_{\kappa(\bar{\nu})}=\name{d}_{\kappa(\bar{\nu})}$, and we let $f(\bar{\nu})=\name{d}_\kappa(\bar{\nu})$ for each $\bar{\nu}\in A^*=\Delta\{A_{\bar{\nu}}:\bar{\nu}\in A^p\}$.
We define $r=p_0^\frown(\bar{U},A^*)$ and conclude that $p\leq^* r$.
Notice that $r+\langle\bar{\nu}\rangle\Vdash\name{c}_{\kappa(\bar{\nu})}=f(\bar{\nu})$ for every $\bar{\nu}\in A^*$.
Applying $\jmath$ we can see that:
$$
\forall\tau\in\kappa^+, \jmath(r)+\langle\bar{U}\upharpoonright\tau\rangle\Vdash \jmath(\name{c})_\kappa=\jmath f(\bar{U}\upharpoonright\tau).
$$
Observe that $\jmath(\name{c})_\kappa$ is a name of a coloring from ${}^{<\kappa}2$ into $2$.
For every $\tau\in\kappa^+$ there is an antichain $\mathcal{A}_\tau$ which describes the possible interpretation of $\jmath f(\bar{U}\upharpoonright\tau)$ as a coloring from ${}^{<\kappa}2$ into $2$.
Since $|\mathcal{A}_\tau|\leq\kappa$ we have a function $g_\tau\in{}^\kappa 2$ in $V$ which guesses simultaneously all the interpretations of $\jmath(\name{c})_\kappa$ mentioned in $\mathcal{A}_\tau$, due to Lemma \ref{lemsimul}.

Let $E=\{g_\tau:\tau\in\kappa^+\}$.
It follows that $\jmath(r)$ forces that $E$ is a witness to $\kappa^+-\Phi_\kappa$ in $M$.
By elementarity, there is a condition $q\geq p$ in $V$ which forces the same statement, so we are done.

\hfill \qedref{thmq1.7}

We indicate that the above weak version of the weak diamond holds in Woodin's model even if $2^\kappa=2^{\kappa^+}$.
We believe that the assumption $\kappa^+<2^\kappa<2^{\kappa^+}$ yields $\neg\Diamond_\kappa\wedge\Phi_\kappa$ in Woodin's model, thus giving a full answer to Question \ref{qdwd}.
However, we still do not know how to prove it.

On the other hand, since our only way to kill diamond at strongly inaccessible cardinals is based on Woodin's method, we always remain with $\lambda^+-\Phi_\lambda$ in the generic extension.
It would be interesting to find out whether this weak form of the weak diamond holds, in \textsf{ZFC}, at every strongly inaccessible cardinal.

\newpage 

\bibliographystyle{amsplain}
\bibliography{arlist}

\end{document}